\documentclass[11pt]{amsart}

\usepackage{amsfonts}
\usepackage{amssymb}
\usepackage{amsmath, amsthm, latexsym}
\usepackage[all]{xy}
\usepackage{hyperref}
\usepackage[margin=1.5in]{geometry}
\usepackage{color}
\def\ms{\medskip}
\def\ni{\noindent}
\def\de{\delta}

\def \c{\mathbb{C}}
\def \z{\mathbb{Z}}
\def \r{\mathbb{R}}
\def \p{\varphi}

\def \T{\mathcal{T}}
\def\B{\mathcal B}
\def\L{\mathcal L}
\def\N{\mathcal N}
\def\K{\mathcal K}
\def\W{\mathcal W}
\def \M{\mathfrak{M}}
\def \b{\mathfrak{b}}
\def \h{\mathfrak{h}}

\def \t{\mathfrak{t}}
\def \g{\mathfrak{g}}
\def\s{\mathfrak{s}}
\def\m{\mathfrak{l}}

\def\R{\c[\s]}
\def\A{\alpha}
\def\st{such that }

\def\iff{if and only if }
\def\gr{\text{Gr}}
\def\coh{cohomology}
\def\eq{equivariant}
\def\ec{equivariant cohomology}
\def\crm{cohomology restriction map}

\def\ef{equivariantly formal}

\def\th{\Theta}
\def\ms{\medskip}
\def\ni{\noindent}

\theoremstyle{plain}
\newtheorem{Th}{Theorem}[section]
\newtheorem{Lem}[Th]{Lemma}
\newtheorem{Prop}[Th]{Proposition}
\newtheorem{Cor}[Th]{Corollary}

\theoremstyle{definition}
\newtheorem{Ex}[Th]{Example}
\newtheorem{Def}[Th]{Definition}
\newtheorem{Rem}[Th]{Remark}

\begin{document}

\title{
torus actions, localization and induced representations on cohomology}
\author{Jim Carrell}
\address{Department of Mathematics, University of British Columbia, Vancouver, B.C., Canada}
\email{jbcarrell@gmail.com}

\noindent
\date{\today}

\begin{abstract} 
This note is  motivated by the problem of understanding 
Springer's remarkable action of the Weyl group $W=N_G(T)/T$ of a semisimple complex
linear algebraic group $G$, with maximal torus $T$, on the \coh\ algebra of an arbitrary
Springer variety in the flag variety of $G$ from the viewpoint of torus actions.
Continuing the work \cite{CK}
which gave a sufficient condition for a group $\W$ acting on the fixed point set of
an algebraic torus action $(S,X)$ on a complex projective 
variety $X$ to lift to a representation of $\W$
on the \coh\ algebra $H^*(X)$ (over $\c$), we
describe when the representation
 on $H^*(X)$ is equivalent to the representation of $\W$ on the \coh\ $H^*(X^S)$
of the fixed point set. As a  consequence of this theorem, we give a simple proof in type $A$ of
the Alvis-Lusztig-Treumann Theorem, which describes Springer's representation of
$W$ for Springer varieties corresponding to nilpotents in a Levi subalgebra of Lie$(G)$.
In the final two sections, we describe the  local structure of the moment graph $\M(X)$
of a special torus action $(S,X)$,  and we also show that if a finite group $\W$ acts on the
moment graph of $X$,  then $\W$ induces pair of actions on $H^*(X)$,
namely the left and right or dot and star  actions of Knutson \cite{knu} and Tymoczko \cite{tym} respectively.
In particular,  $W$ acts on the moment (or Bruhat) graph $\M(G/P)$ of $(T,G/P)$ for any parabolic $P$ in $G$ containing $T$,
and the right action of $W$ on $H^*(G/P)$ is an induced representation. 
Furthermore, we show the left action
of $W$ on $H^*(G/P)$ is trivial.
\end{abstract}

\maketitle

\section{introduction}

In \cite{Spr2,Spr3}, T. A. Springer constructed a remarkable representation of the Weyl group $W$ of a semisimple complex linear algebraic group $G$ on the \coh\ algebra $H^*(X)$ 
with coefficients in  $\c$ of an arbitrary
Springer variety $X$ in the flag variety $G/B$. One remarkable aspect of Springer's 
representation is that $W$ does not act on $X$ itself, except when 
$X=G/B$. Springer's original construction uses $\ell$-adic \coh, but, later,  
more direct constructions for ordinary cohomology 
became available  (cf. \cite{bbm,COMP,KL,slod,Spr4}). 
In \cite{CK} (also see \cite{COMP}), the author and K. Kaveh showed that 
Springer's action when $G=SL(n,\c)$ is a consequence of a general result on torus actions $(S,X)$, where
$X$ is a projective variety with no odd \coh. This result gives a sufficient condition for the action of a finite group 
$\W$ on the \coh\ algebra of the fixed point set $X^S$ of $S$ to lift to an action on the \coh\ algebra of $X$.

The problem we consider here is the natural question of determining which representations 
of $\W$ can occur on $H^*(X)$ as a result of lifting a representation on $H^*(X^S)$.
In particular, our result describes when the representations of $\W$ on $H^*(X^S)$ and $H^*(X)$ are in fact equivalent. 
For example,  if the representation of $\W$ on $H^*(X^S)$ comes from an action of $\W$
on $X^S$ itself, then the representation of $\W$ on $H^*(X)$ is determined by the $\W$-orbits 
on $X^S$. If the action $(\W,X^S)$ is simply transitive, then the 
representation of $\W$ on $H^*(X)$ has to be the regular representation. 
Applying this result to Springer's Weyl group representation in the case of $SL(n,\c)/B,$
and more generally to Springer varieties $X$ of parabolic type for 
which the \crm\ associated to the inclusion
$X\to G/B$ gives a surjection $H^*(G/B)\to H^*(X)$, we get a 
short proof elementary proof of the Alvis-Lusztig-Treumann Theorem \cite{AL,Treu}.

\section{preliminary remarks and statement of the main result}
Throughout this note, $X$ will denote a connected complex projective variety with 
vanishing odd \coh\ on which an algebraic torus $S=(\c^*)^k$ acts algebraically. In particular,
the action $(S,X)$ is \ef.
We remark that many of the results in this note are also valid in the topological setting for 
\ef\ $(S^1)^k$-actions. However, some of the results in the projective case
are sharper, and  the applications 
involve complex projective algebraic varieties. Thus we will work in that setting.

Let us briefy recall some definitions. If $Y$ is a 
topological space, then $H^*(Y)$ will  denote the \coh\
algebra of $Y$ with coefficients in $\c$. Assume $Y$ admits a 
topological action by the algebraic torus $S$.
The $S$-\ec\ of $Y$ over $\c$ is defined to be
$H^*_S(Y)=H^*(Y_S)$, where $Y_S$ is the Borel space $(Y\times E)/S$, 
$E$ being a contractible space with a free action of $S$. 
Note that $H^*_S(Y^S)=H^*(E/S)\otimes H^*(Y^S)$, by the Kunneth Formula.
Since the \coh\ algebra  $H^*(E/S)$ can be  identified via
the characteristic homomorphism with the graded algebra
$\c[\s]$, where $\s=\text{Lie}(S)$, we will make the identification
$H^*_S(Y^S)=\c[\s]\otimes H^*(Y^S)$, where the algebra $\c[\s]$
is assumed to be graded in even degrees. 
The forgetful map $H_S^*(Y)\to H^*(Y)$ is the \crm\
induced by the inclusion of $Y$ into $Y_S$ along a fibre. 
If the action$(S,Y)$ is \ef, one has an exact sequence
\begin{equation}\label{exact-seq}
0\longrightarrow \c[\s]^+ H_S^*(Y) \longrightarrow H_S^*(Y) \longrightarrow H^*(Y) \longrightarrow 0.
\end{equation}
Here $\c[\s]^+$ is the augmentation ideal of $\c[\s]$ (i.e. the maximal ideal of $0\in \s$).

Recall the  Localization Theorem of Atiyah and Bott \cite{AB} (see \cite{BRION}
for a nice proof): 
if the action $(S,Y)$ is \ef, then the inclusion $i:Y^S \to Y$ induces 
an injection $i_S^*:H_S^*(Y)\to H_S^*(Y^S)$, and, furthermore, $i_S^*$ 
becomes an isomorphism after 
localizing $H_S^*(Y^S)$ at finitely many elements  $f_1, \ldots, f_r \in \c[\s]$.

Now suppose $Y$ is  a complex projective variety with vanishing odd \coh,
and suppose the algebraic torus $S$ acts on $Y$.
Let $\W$ be a  finite group acting topologically on $Y$, and 
suppose $\W$ commutes with the action of $S$:
that is, $\W \times S$ acts on $Y$. Then $\W$ also 
acts on $Y^S$, so it acts on $H^*_S(Y)$ and $H^*_S(Y^S)$ 
by $\c[\s]$-module algebra isomorphisms. Suppose
$X$ is an $S$-stable subvariety of $Y$, and note that $X^S$ 
is  nonempty by the Borel Fixed Point Theorem. 
We do not assume, however, that $\W$ acts on either $X$ or $X^S$. 
The main result in  \cite{CK} is the following:

\begin{Th} \label{CK}
Suppose the \crm\ $\iota:H^*(Y)\to H^*(X)$ is surjective, and assume
also that $H^*(X^S)$ is a $\W$-algebra so that the 
\crm\ $j:H^*(Y^S)\to H^*(X^S)$ is $\W$-\eq. 
Then there exists an action of $\W$ on $H^*_S(X)$ via 
$\c[\s]$-module isomorphisms such that 
$\iota_S:H_S^*(Y)\to H_S^*(X)$ is $\W$-\eq, and the 
\crm\ $i_S^*:H^*_S(X)\to H^*_S(X^S)$ is also $\W$-\eq\
for the action of $\W$ on  $H^*_S(X^S)$ via $\c[\s]$-module isomorphisms
induced by the action of $\W$  on $H^*(X^S)$.
In particular, the forgetful map induces a graded $\c$-algebra representation of $\W$ on $H^*(X)$. 
\end{Th}

\ms
Although the \crm\ $j:H^*(Y^S)\to H^*(X^S)$ is not assumed to be surjective, it turns out,
as we will prove in Theorem \ref{thmSURJ}, that surjectivity 
of $j$ follows automatically from the surjectivity of 
$\iota$. The relevance of Theorem \ref{CK} is that it gives a sufficient condition
for an action of $\W$ on $H^*(X^S)$ to lift to $H_S^*(X)$.

We now state our main result. Here, the action $(S,Y)$ isn't needed,
but $X$ is, as usual,   assumed to have vanishing odd \coh.
\begin{Th} \label{2.2}  Suppose a finite
group $\W$ acts on the  $\c$-algebras $H_S^*(X)$ and $H_S^*(X^S)$
via graded $\c[\s]$-module algebra isomorphisms and that 
the \crm\ $H_S^*(X)\to H_S^*(X^S)$ is $\W$-\eq. Thus $\W$ acts  
as graded $\c$-algebra isomorphisms on both $H^*(X)$ and $H^*(X^S)$ via the forgetful maps. 
Moreover, $H^*(X^S)$ has a $\W$-stable filtration 
\begin{equation}\label{eq2}
F_{-1} =\{0\}  \subset F_0 \subset F_1 \subset \cdots \subset F_m 
\subset F_{m+1} \subset \cdots  \subset H^*(X^S)
\end{equation} 
such that 
\begin{equation}\label{eq3}
F_iF_j\subset F_{i+j} ~~\makebox{and} ~~
\gr_F H^*(X^S)\cong H^*(X)
\end{equation} via a $\W$-\eq\  isomorphism of graded $\c$-algebras.
In particular,
$H^*(X)\cong H^*(X^S)$ as $\W$-modules.

\end{Th}

In the case  where $X^S$ is finite, one  gets 
an explicit description of the action of $\W$ on $H^*(X)$.

\begin{Cor} \label{COR1} Assume $X^S$ is finite, and 
the action of $\W$ on $H_S^*(X^S)$ arises from an action of $\W$ on $X^S$.
Then:

\ms
\ni
$(i)$  if  the action $(\W,X^S)$ is simply transitive,  the representation of $\W$
on $H^*(X)$ is the regular representation, and

\ms
\ni
$(ii)$ if the action $(\W,X^S)$ is transitive, 
then the representation of $\W$ on $H^*(X)$  is the induced
representation $\mathrm{Ind}_{\W_z}^{\W}(\c)$ associated to 
the trivial one dimensional representation 
of the isotropy group  $\W_z$ 
for any $z\in X^S$. 

\ms
\ni
In general, the representation of $\W$ on $H^*(X)$  is either the regular
representation or a direct sum of induced representations corresponding to  the orbits of $\W$
on $X^S$.
\end{Cor}

{\color{black} 

\section{The Evaluation mapping} 
As usual, we  suppose  that $X$ has an algebraic torus action $(S,X)$
and vanishing odd \coh. Our plan is to 
obtain an increasing filtration $F_{\cdot}$ of $H^*(X^S)$ with $F_iF_j \subset F_{i+j}$
whose associated graded algebra 
is isomorphic with $H^*(X)$. 

\begin{Th} \label{FILTRATION} Assume $X$ has vanishing odd \coh\ and
let $(S,X)$ be an algebraic torus action.  Then $H^*(X^S)$ admits a filtration 
\begin{equation} \label{filtn}
F_{-1}=\{0\} \subset F_0 \subset F_1 \subset \cdots  
\subset F_m \subset F_{m+1} \subset \cdots \subset H^*(X^S)
\end{equation}
satisfying $F_iF_j \subset F_{i+j}$ with the property that $F_i/F_{i-1} \cong H^{2i}(X)$
for all $i\ge 0$, and 
the following map is a graded $\c$-algebra isomorphism: 
\begin{equation} \label{eqCOHOM} 
\mathrm{Gr}_F H^*(X^S) =\bigoplus _{i\ge 0} F_i/F_{i-1} \cong \bigoplus_{i\ge 0} H^{2i}(X)=H^*(X).
\end{equation}
\end{Th}

\proof
Since  the cohomology $H^*(X)$ is trivial in odd degrees, it follows
that the action of $S$ is \ef. Thus
$H_S^*(X)$ is a free $\R$-module of rank $r=\dim H^0(X^S)$.
Let $f_1, \ldots, f_r \in \R$ be elements \st the restriction
 $i_S^*:H_S^*(X)\to H_S^*(X^S)$ is an isomorphism after 
localization at these elements.
We will say $a\in \s$ is {\em regular} if $f_i(a)\ne 0$ for $i=1, \dots r$.
By equivariant formality, we have the following well known exact sequence
\begin{equation}\label{exact-seq}
0\longrightarrow \R^+ H_S^*(X) \longrightarrow H_S^*(X) \longrightarrow H^*(X) \longrightarrow 0,
\end{equation}
where the map $H_S^*(X) \longrightarrow H^*(X)$ is induced by inclusion of $X$ in $X_S$
and $\R^+$ is the augmentation ideal (i.e. the maximal ideal of $0\in \s$).
If $V$ is a $\c$-vector space and $a\in \s$, let $V_a$ be the $\R$-module 
defined by putting $f\cdot v=f(a)v$
for any $f\in \R$ and $v\in V$. 
Similarly, if $M$ is an $\R$-module, put $M[a]=M\otimes_{\R} \c_a$.
When $M$ is free of rank $r$, then $\dim _\c M[a]=r.$ 
  Now, the exact sequence (\ref{exact-seq}) implies
$$H_S^*(X)[0]\cong H^*(X)[0]=H^*(X).$$ 
Furthermore, the Localization Theorem implies that 
\begin{equation}\label{isom}
i_S^*: H_S^*(X)[a]\to H_S^*(X^S)[a]
\end{equation}
is an isomorphism of $\R$-algebras for any regular $a\in \s$, hence also an isomorphism of 
(ungraded) $\c$-algebras. 
For any $i\ge 0$, put 
$${\mathcal F}_i=\sum _{j\le i} H_S^{2j}(X).$$
This defines an increasing filtration of $H_S^*(X)$ such that 
 ${\mathcal F}_i {\mathcal F}_j\subset {\mathcal F}_{i+j}$ for $i,j\ge 0$.
Moreover, defining ${\mathcal F}_i[a]$ to be the image of ${\mathcal F}_i$
in $H_S^*(X)[a]$, we obtain, for any regular $a$, a filtration 
$F(i,a)$ of $H^*(X^S)$ as follows.  Let $e_a: H_S^*(X^S)[a] \to H^*(X^S)$ 
be the unique $\c$-algebra isomorphism such that 
$e_a(f\otimes \phi)=f(a)\phi$. Now put
$$F(i,a)=e_a(i_S^*({\mathcal F}_i[a])).$$
Then define the filtration (\ref{filtn}) of $H^*(X^S)$ to be this filtration.
Returning to the filtration ${\mathcal F}_i[a]$ of $H_S^*(X)[a]$,
we obtain from the forgetful map in (\ref{exact-seq}) 
a filtration $G(i,a)\subset G(i+1,a)$
of $H^*(X)=H^*(X)[a]$ such that   
$$\bigoplus _{i\ge 0} G(i,a)/G(i-1,a) \cong \bigoplus _{i\ge 0} F(i,a)/F(i-1,a)$$ 
is a graded $\c$-algebra isomorphism via (\ref{isom}).
It remains to establish the graded ring isomorphism
$$\bigoplus_{i\ge 0} G(i,a)/G(i-1,a) \cong H^*(X).$$
But $G(i,a)/G(i-1,a)\cong H^{2i}(X)[a]$, and 
this gives rise to  the natural multiplication 
$$G(i,a)/G(i-1,a)\otimes G(j,a)/G(j-1,a)\to H^{2i}(X)[a]\otimes H^{2j}(X)[a] \to H^{2(i+j)}(X)[a]$$
via  the ring structure of $H^*(X)[a]$. But $H^*(X)=H^*(X)[a]$ by definition, 
so the proof of Theorem \ref{FILTRATION} is finished. \qed

\medskip

\begin{Rem} The filtration $F_{\cdot}$ is not related to the natural grading of $H^*(X^S)$.
\end{Rem}

\begin{Rem} The fact that 
$H^*(X^S)$ has a filtration whose associated graded is $H^*(X)$ is also 
shown by Puppe (cf. \cite[p.14]{VP2}) for $S^1$-actions on compact spaces.
Previously, in \cite{CL}, Lieberman and the author had proven this result
smooth projective varieties having a holomorphic vector field with simple zeros,
which includes the case of an arbitrary algebraic torus action on $X$.
In the compact Kaehler case,  (5) was sharpened in \cite{CKP} by bringing in the  Hodge 
decomposition of the \ec\ $H^*_S(X)$.
\end{Rem}

Since the filtration (6) is preserved by equivariant maps,  one gets 
the following interesting consequence which we believe may be new. 

\begin{Th} \label{thmSURJ}  Assume  that $X$ and $Y$ are projective with 
vanishing odd \coh\ and algebraic torus actions $(S,X)$ and $(S,Y)$,
and suppose $f:X\to Y$ be a regular $S$-\eq\ map such that 
$f^*:H^*(Y)\to H^*(X)$ is surjective. Let $\rho$ denote the
restriction of $f$ to $X^S$.
Then $\rho^*: H^*(Y^S) \to H^*(X^S)$ is also surjective.
\end{Th}

\proof Let $E_\cdot$ and $F_\cdot$ be the filtrations of $H^*(Y^S)$ and $H^*(X^S)$ 
defined in Theorem \ref{FILTRATION}. Then    
$\rho^*(E_i)\subset F_i$ for all $i$. This gives the commutative diagram

\begin{equation} \label{equ-comm-diag}
\xymatrix{
\gr_E H^*(Y^S) \ar[r] \ar[d]^{\hat{\rho}} & H^*(Y) \ar[d] ^{f^*}\\
\gr_F H^*(X^S) \ar[r] & H^*(X),\\
}
\end{equation}
where $\hat{\rho}$ is the graded $\c$-algebra morphism induced by $\rho^*$. 
Since the horizontal maps are isomorphisms
and $f^*$ is surjective, it follows that $\hat{\rho}$ is surjective, so 
$\rho^*(E_i)=F_i$ for all $i$.
Since $\bigcup F_i =H^*(Y^S)$, it follows that $\rho^*$ is surjective. \qed

\begin{Rem} Of course, the above hypotheses  immediately imply that 
$\chi(Y)\geq \chi(X)$ for the Euler characteristics of $X$ and $Y$. The theorem 
strengthens this fact.
We will  apply it in the proof of the 
Alvis-Lusztig-Treumann Theorem in Section \ref{secSPRINGER}.
\end{Rem}

\section{the Proof of Theorem \ref{2.2}}
Note first that, by assumption, the map $i_S^*: H_S^*(X)[a]\to H_S^*(X^S)[a]$
is $\W$-\eq\ with respect to the action of $\W$ on $H_S^*(X^S)[a]$ 
for any regular $a\in \s$. In addition, the map 
$e_a: H_S^*(X^S)[a] \to H^*(X^S)$ is also $\W$-\eq.
Thus, putting  $F_i=F(i, a)$ gives a  $\W$-\eq\ filtration of $H^*(X^S)$
having the required properties whose associated graded is $H^*(X)$
by Theorem \ref{FILTRATION}. 
It remains to show that the representation of $\W$
on $H^*(X^S)$ is equivalent to the representation of $\W$ on  
$\mathrm{Gr}_F H^*(X^S)=\bigoplus _{i\ge 0} F_i/F_{i-1}$. This follows from the following lemma.
\begin{Lem} Let $\W$ be a finite group acting linearly on a $\c$-vector space 
$V$ having an $\W$-invariant filtration 
$$V_0 =\{0\} \subset V_1 \subset \cdots V_{m-1} \subset V_m =V.$$
Then the induced representation of $\W$ on $Gr~V=\bigoplus _{i\ge 0}V_i/V_{i-1}$ 
is equivalent  to the given representation of $\W$ on $V$.
\end{Lem}
\proof Since $\W$ is finite, every $\W$ invariant subspace of $V$ has an $\W$-invariant 
complement. Applying this fact to $F_i \subset F_{i+1}$ for each $i$, one gets a $\c$-linear
$\W$-\eq\ isomorphism between $V$ and  $Gr~V$. \qed

\medskip
To prove Corollary \ref{COR1}, suppose $X^S=\{x_1, \dots, x_r\}$ and 
$\W$ acts transitively. Let $\K$ be the isotropy group of $x_1$. 
By the Orbit Stabilizer Theorem, $|\W/\K|=r$.
Since $H^0(X^S)=\c^{X^S}$, a basis of $H^0(X^S)$ is 
given by the functions $\de_1, \dots, \de_r$ on $X^S$ defined by
the conditions $\de_i(x_j)=\de_{ij}$. By definition, 
the subspace $Y=\c \de_1$ is a $\K$-stable line. Choose $w_1, \dots, w_r \in \W$
so that $w_i(x_1)= x_i$. Then $\W/\K=\{w_1\K, \dots, w_r\K\}$, 
and $w_i \cdot \de_1= \de_i$. For $w_i \cdot \de_1 (x_i)=\de_1((w_i)^{-1} \cdot x_i)=\de_1(x_1)=1$
while if $i\neq j$, $(w_i)^{-1} \cdot x_j \ne x_1$, so
$w_i \cdot \de_1 (x_j)=0$. Consequently, 
$$H^0(X^S)=\sum_{\sigma \in \W/\K} \sigma Y.$$
If $\K=\{e\}$, then $(\W, H^0(X^S))$ is the regular representation.
Otherwise, the representation of $\W$ on $H^0(X^S)$ is induced as claimed. 
The remaining assertion is straightforward. \qed

\section{Springer's representation of $W$ on a Springer variety}\label{secSPRINGER}

In this section, we will consider Springer's Weyl group action 
on the \coh\ of a Springer variety. We will begin with a brief introduction to 
the flag variety $G/B$ followed by some remarks about Springer varieties.
Let $G$ be a semi-simple complex algebraic group, a suppose
$B$ and $T$ denote, respectively, a Borel subgroup of $G$ and a
maximal torus of $G$ contained in $B$. Also, let $\g$, $\b$ and $\t$ 
denote their corresponding Lie algebras, and let $W= N_G(T)/T$ of $(G,T)$
be the Weyl group of the pair $(G,T)$.
The flag variety of $G$ is the homogeneous space  $\B=G/B$.
It is a complex projective variety which plays 
an important role in representation theory due to the fact that it
parameterizes the set of all Borel subgroups of $G$ via the fundamental bijection
$gB\to gBg^{-1}$. Equivalently, $\B$ parameterizes
the set of all Borel subalgebras of $\g=T_{B}(\B)$ via $gB\to \text{Ad}(g) \b$. 
The Bruhat decomposition of $G$, namely $G=\cup_{w\in W} BwB$ induces a CW decomposition
of $\B$ with even dimensional cells each isomorphic to a complex affine space, 
so $\B$ has vanishing odd \coh. Moreover, $T$ defines an algebraic torus action 
on $\B$ by left translation. Since
$\B$ can be identified topologically with $K/H$, where $K$ is a maximal compact subgroup of $G$
and $H=K\cap T$, and since $W\cong N_K(H)/H$, it follows that $W$ acts topologically on $\B$ via 
$$w\cdot kH=k \dot{w}^{-1}H,$$ 
where $\dot{w}\in N_K(H)$ is a representative of $w$.
This action commutes with the action of $T$ on $\B$,
and allows $\B^T$ to be naturally identified with $W$ via $w\to \dot{w}B$.

Let $\N$ denote the cone of all nilpotent elements in $\g$. 
The {\em Springer variety associated to} $x\in \N$
is the subvariety $\B_x$ of $\B$ consisting 
of all Borel subalgebras of $\g$ containing $x$; 
$\B_x$  is connected but not irreducible except when $x=0$ or $x$ is 
a principal (or regular) nilpotent: that is, $C_G(x)$ has minimum dimension. 
The Springer variety $\B_x$ may also be defined as the fixed point set $\B^u$ in $\B$ of
 the unipotent element $u=\text{exp}(x)$ which, by the above identification, is the set 
 of all Borel subgroups of $G$ containing $u$. A famous result of DeConcini-Lusztig-Procesi \cite{DLP}
 is that $H^{odd}(\B_x)$ is trivial for all $x\in \N$. Moreover, it is well known
that the fixed point set of a unipotent group acting on a projective variety is always conected.
Furthermore, every Springer variety $\B_x$ in $\B$ when $G=SL(n,\c)$
admits an affine paving which is part of an affine paving of $\B$ \cite{SPALT}, so
the \crm\ $H^*(\B)\to H^*(\B_x)$ is surjective for all $x\in \N$. 

A celebrated result of Springer (\cite{Spr2,Spr3}) says that 
for every $x\in \N$, there exists 
a graded $\c$-algebra representation of $W$ on $H^*(\B_x)$.
Since $W$ doesn't act on $\B_x$ itself, except in a few special cases,
his original construction uses a good deal of heavy machimery.
There are now much more elementary approaches, using, for example,  intersection \coh\ or homotopy theory
\cite{bbm,KL,slod,Treu}. A natural question, of course, is how to explicitly classify 
these representations. Every Springer variety $\B_x$ admits a torus action (by the 
Jacobson-Morosov Lemma), but it isn't known in general when $W$ acts on $(\B_x)^S$,
so Theorem 2.2 doesn't help.
However,  Proposition \ref{ALT} below answers this when $H^*(\B)\to H^*(\B_x)$ is surjective.
On the other hand, the Alvis-Lusztig-Treumann Theorem classifies Springer's 
representations for nilpotents in a Levi subalgebra of $\g$. 
Consider a Levi subgroup $L=C_G(S)$ of $G$, where $S$ is a subtorus of $T$,
and suppose  $\m$  denotes the Lie algebra  of $L$. Also, let $W_L$ denote the Weyl group of $(L,T)$.
The following lemma reviews some basic facts 
about the action $(S,\B_x)$.when $x\in \N \cap \m$.

\begin{Lem} Let $S$ be a subtorus of $T$ and $L=C_G(S)$. Then:

\ms
\ni
$(i)$ every  component of the fixed point set $\B^S$ may be identified with a flag variety of the form 
$L/(L\cap wBw^{-1} )$ for some $w\in W$;

\ms
\ni
$(ii)$ in   fact,   $\B^S$ has $|W/W_L|$ components;

\ms
\ni
$(iii)$ if  $x\in \N\cap \m$ is a principal nilpotent in $\m$, then 
every component of $\B^S$ contains a unique point if $(\B_x)^S$, 
hence $|(\B_x)^S|=|W/W_L|$; and

\ms
\ni
$(iv)$ if $x$ is as in $(iii )$, then $W$ acts transitively on $(\B_x)^S$ itself.
\end{Lem}
\proof $(i)$ follows from the bijection $W\to \B^T$ defined above. For $(ii)$, see 
\cite[Lemma 6.3]{JBC}. For $(iii)$ and $(iv)$, use the proof of Theorem 2 in \cite{COMP}.
\qed

\ms Note that we will view each flag variety $L/L\cap wBw^{-1}$ as a subvariety
of $\B$, namely a comonent of $\B^S$.
The next Proposition, which seems to be of some independent interest,  will be helpful in the proof of the 
Alvis-Lusztig-Treumann Theorem. As usual, $L=C_G(S)$.
\begin{Prop}\label{ALT}  Let $x\in \N \cap \m$ and assume  the \crm\ $H^*(\B) \to H^*(\B_x)$ is 
surjective. If $x$ is a principal nilpotent in $\m$, then, 
for any subtorus of $S'$ of $S$, $W$ acts on $H^*((\B_x)^{S'})$ 
so that the \crm\ $H^*((\B)^{S'}) \to H^*((\B_x)^{S'})$ is $W$-\eq.
\end{Prop}

\proof Since $x$ is a principal nilpotent in $\m$, it follows that
$W$ acts on $(\B_x)^S$ from part $(iv)$ of the previous lemma.
Let $S'$ be any subtorus of $S$. By assumption,
the \crm\ $H^*(\B) \to H^*(\B_x)$ is 
surjective, hence by Theorem \ref{thmSURJ}, the
\crm\ $H^*((\B)^{S'}) \to H^*((\B_x)^{S'})$ is also surjective.
 But $S$ gives a torus action on both $(\B)^{S'}$ and $(\B_x)^{S'}$,
 so $W$ acts on $H^*((\B_x)^{S'})$ as claimed by Theorem \ref{CK}. \qed
 
 \ms
Put $\B_L=L/L\cap B$, and,
for any nilpotent $x\in \m$, let $\B_{L,x}$ denote 
the Springer fibre associated to $x$ in the flag variety $\B_L$.
Then the Alvis-Lusztig-Treumann Theorem  \cite{AL,L,Treu} (over $\c$) goes as follows:
 \begin{Th} \label{THMAL} Suppose $x\in \N \cap \m$. Then the $W$-action 
 on $H^*(\B_x)$ is equivalent to the induced representation 
 $\mathrm{Ind}_{W_L}^W(H^*(\B_{L,x}))$.
 \end{Th}
 
In the following (partial)  proof, we will assume  the \crm\ $H^*(\B) \to H^*(\B_x)$ is surjective
and $x$ is principal in a Levi subalgebra $\tilde{\m}$ which is the Lie algebra of $\tilde{L}=C_G(\tilde{S})$.
Thus, $(\B_x)^{\tilde{S}}$ is finite. (As pointed out below, this assumption holds in type $A$.)   
Now consider a subtorus $S\subset \tilde{S}$ with $L=C_G(S)$  such that $x\in \m$. 
The Weyl group $W_L=N_L(T)/T$ of $L$ acts on the flag variety 
 $\L_w=L/(L\cap wBw^{-1})$ for each $w\in W$. 
Letting $w_1, \dots, w_k\in W$ denote a complete set of
 representatives for $W/W_L$ and putting $\L_{w,x}=\B_x \cap (L/L\cap wBw^{-1})$, we have
 \begin{equation}\label{9}
 (\B_x)^S=\bigcup_{1\le i \le k} \L_{w_i,x}.
 \end{equation}
 By  Proposition \ref{ALT}, $W$ acts on $H^*((\B_x)^S)$, and
 by Theorem \ref{2.2} there is a $W$-stable filtration $F_{\cdot}$
 of $H^*((\B_x)^S)$ such that $\gr _F H^*((\B_x)^S) \cong H^*(\B_x)$
 as graded $W$-algebras.  Since
$$H^*((\B_x)^S)=\bigoplus _{1\le i \le k} H^*(\L_{w_i,x}),$$
it follows that $H^*((\B_x)^S)\cong \mathrm{Ind}_{W_L}^W(H^*(\B_{L,x}))$,
due to the fact that $\L_{e,x}=\B_{L,x}.$
\qed

\ms
By the above remarks, we have a more explicit description in type $A$, 
since, when $\g=\s\m(n,\c)$, every nilpotent $x$ is a principal nilpotent, 
and $H^*(\B)\to H^*(\B_x)$ is always surjective. 
   
\begin{Cor}\label{5.3} Suppose $G=SL(n,\c)$ and $x\in \N$ is a principal nilpotent for $\m$.  
Then Springer's representation of $W$ on $H^*(\B_x)$ is equivalent to 
$\mathrm{Ind}_{W_L}^W (\c)$.  
\end{Cor} 

\begin{Rem} Springer's representation in type $A$
has been related to nilpotent orbits in $\g=\s\m(n,\c)$. 
In fact, if $x$ and $y$ are nilpotents in $\g$ which are dual in
the sense the partitions associated to $x$ and $y$ are dual to each other,
DeConcini and Procesi \cite{DP}  showed that there is a $W$-\eq\
isomorphism of graded $C$-algebras
$H^*(\B_x)\cong A(\overline{O_y} \cap \t)$, where $O_y$ is the 
$G$-orbit of $y$ in $\g$ and $A(\overline{O_y} \cap \t)$ is the coordinate ring of the
schematic intersection of $\overline{O_y}$ and the Cartan subalgebra  $\t$. 
Since the Weyl group, $S_n$ in this case, acts on this coordinate ring
as an induced representation, this gives an alternative to the Alvis-Lusztig-Treumann 
Theorem for type $A$. 
This result verifies an important case of a general conjecture of 
Kraft \cite{kra} for pairs of dual nilpotents. Kraft's conjecture 
is not true in general, but an alternative and more general proof of 
some cases Kraft's conjecture including the $\g=\s\m(n,\c)$ case was
given by the author in \cite{JBC} using torus actions and Weyl group orbits.
\end{Rem}

\medskip
\begin{Rem}  
The general question of when the \crm\ $H^*(\B)\to H^*(\B_x)$ is 
surjective is interesting and unsolved. A necessary condition for surjectivity
in the middle dimension is that the component group $A(x)=C_G(x)/C_G(x)^0$ 
of $G$ is trivial.
The conjugacy classes of these groups have been computed in  \cite{somm}, where
one sees there are in general many non-surjective examples. 
It has recently been  shown that
if $G$ is simply laced, then $H^2(\B)\to H^2(\B_x)$ is always surjective \cite{CVX},
but this is only helpful when $H^2(\B_x)$ generates $H^*(\B_x)$.
\end{Rem} 

\section{Special torus actions and the moment graph}

In this section and the next, we will study a class of algebraic torus actions $(S,X)$ called special in \cite{cp}.
These actions have the property that $X^S$ is finite, and $X$ contains only finitely many $S$-stable curves:
that is, curves which are the closure of a one dimensional $S$-orbit.
The set of these curves will be denoted by $E(X)$. Our plan in this section is to point out some   basic properties of the moment graph $\M(X)$, namely the unoriented graph 
with vertex set $X^S$ such that the set of edges at each vertex $x$ consists of the set 
$E(X,x)=\{C\in E(X) \mid x\in C^S\}$.
Note that for any $x\in X^S$, $S$ acts on the Zariski tangent space
$T_x(X)$  of $X$ at $x$ via the differential of the morphism $z\to s\cdot z$
at $z=x$.

\begin{Def} An  algebraic torus action $(S,X)$ is called {\em special} if and only if:

\ms
\ni
$(i)$ $(S,X)$ is locally linearizable;

\ms
\ni
$(ii)$  for any $x\in X^S$, $\dim T_x(X)^S=0$; and

\ms
\ni
$(iii)$ for any codimension one torus $H\subset S$, $\dim T_x(X)^H\le 1$. 

\end{Def}

\ms
If $(S,X)$ is special, then for every $S$-stable subvariety $Y$ of $X$, $(S,Y)$  is
also special.  The basic examples of special actions we will consider are the actions $(T,G/P)$ where  $T$ is a maximal torus in the Borel $B$ and $P$ is a parabolic subgroup of $G$ containing  $B$ and their $T$-stable subvarieties. We will describe the moment graph of $(T,\B)$ below.  

The next lemma describes some of the
properties of $S$-stable curves when $(S,X)$ is special.
\begin{Lem} Let $(S,X)$ be special. Then $X^S$ is finite, and $X$ contains only finitely many
$S$-stable curves. Moreover, every $S$-stable curve $C$ is smooth and $C^S=\{x,y\}$ for a pair of 
distinct points $x,y \in X^S$. If $C_1, \dots, C_k$ are distinct $S$-stable curves containing $x\in X^S$, then
the sum
$$\sum_{1\le i\le k} T_x(C_i)\subset T_x(X)$$ 
is direct. Finally, if $C$ is an $S$-stable curve in $X$,  $C^S=\{x,y\}$, 
and if the representation of $S$ on $T_x(C)$ is given by the character $\chi_C$,
then the representation of $S$ on $T_y(C)$ is given by $(\chi_C)^{-1}$.
\end{Lem}

\ms
We will omit the proof. An $S$-stable
curve $C$ is rational, hence is homeomorphic to 
$S^2$ with a smooth $(S^1)^k$-action ($k=\dim S$) whose fixed point set  is $C^S$. 
Our  next result gives a picture of the moment graph of a special torus action
at an arbitrary vertex in the sense that  we will give geometric bounds, above and below, 
for the number $|E(X,x)|$ of edges of $\M(X)$ at $x$. We begin with a uniform lower bound. 

\begin{Lem} Given a special action $(S,X)$, we have 
$$|E(X,x)|\ge \dim_x X$$
for every $x\in X^S$. Thus, if $X$ is irreducible, then
$\M(X)$ has at least $\dim X$ edges at each vertex. 
\end{Lem}

\proof This is proved in \cite[Lemma, P. 56]{cp}. \qed

\medskip
We will next obtain an upper bound for $|E(X,x)|$, which only depends on the geometry of $X$.
Let $X(S)$ and $Y(S)$ denote respectively the group of characters of $S$ 
and the group of one parameter
subgroups of $S$. Recall the canonical perfect pairing $X(S)\times Y(S) \to \z$ defined by 
$<\A,\lambda>=m$ \iff $\A(\lambda(t))= t^m$.
For this bound, we need to consider the reduced tangent cone $\T_x(X)\subset T_x(X)$ 
of $X$ at $x$ and its  linear span $\th_x(X)$
in $T_x(X)$. Also let $\theta_x(X)\subset X(S)$
denote the weights of $\th_x(X)$.  
The structure of the moment graph at an arbitrary vertex is given by the following:

\begin{Th} Assume  $(S,X)$ is special. Then,  for any $x\in X^S$,  
\begin{equation}\label{10} 
 \dim_x X \le |E(X,x)|\le \dim \th_x(X). \end{equation}
Moreover:

\medskip
\noindent
$(i)$ if $X$ is smooth at $x$, then both inequalities in 
$(10)$ are equalities;

\medskip
\noindent
$(ii)$ if all weights of $\th_x(X)$ are extreme points of their convex hull
in $X(S)\otimes \r$, then every $S$-stable line in $\th_x(X)$ is contained in $\T_x(X)$; and

\medskip
\noindent
$(iii)$ if every $S$-stable line in $\T_x(X)$ is $T_x(C)$
for some $C\in E(X,x)$ and the convexity condition on the $S$-weights 
on $\th_x(X)$ in $(ii)$ holds, then $|E(X,x)|=\dim \th_x(X)$.
\end{Th}

\medskip
\proof Let $x$ be a smooth point. Then $\dim X_x \le |E(X,x)|\le \dim T_x(X)=\dim_x X,$
since the sum $\sum_{C\in E(X,x)} T_x(C)\subset T_x(X)$ is direct. Thus $(i)$ follows.
For $(ii)$, let $L_\A$ be an $S$-stable line of weight $\A$ in $\th_x(X)$, and let $y\in L_\A$
be nonzero. We may write write $y=z_1+ \cdots +z_r$ where each $z_i \in \T_x(X)$. Since each $z_i$ 
lies in  $\T_x(X)$,  we may decompose each $z_i$ into a sum of $S$-weight vectors
in $\th_x(X)$, say $z_i=\sum_{\beta \in \theta_x(X)} w_{i,\beta}$.
Since $\A$ is an extreme point of the convex hull of $\theta_x(X)$ in $X(S)\otimes \r$, 
there exists a $\lambda \in Y(S)$ such that $<\beta -\A,\lambda>>0$ for all $\beta \in \theta_x(X)$
such that $\beta \ne \A$. Now
\begin{eqnarray*}
t^{-<\A,\lambda> }\lambda(t)\cdot z_i&=& t^{-<\A,\lambda> } \sum_{\beta} \lambda(t)\cdot w_{i,\beta} \\
&=&t^{-<\A,\lambda> } \sum_{\beta} t^{<\beta,\lambda>} w_{i,\beta} \\
&=&w_{i,\A} + \sum_{\beta \ne \A} t^{<\beta -\A, \lambda>} w_{i,\beta}
\end{eqnarray*}
Hence, $\lim_{t\to 0}t^{-<\A,\lambda> }\lambda(t)\cdot z_i=w_{i,\A}$. 
It follows that $y=\sum_{i} w_{i,\A}$, and, since $y\ne 0$, we conclude
some $w_{i,\A}\ne 0.$ But since 
$\T_x(X)$ is an $S$-stable closed cone,  each $w_{i,\A}$ lies in 
$\T_x(X)$. Hence  $L_\A\subset \T_x(X)$. Part $(iii)$ follows immediately from $(ii)$. \qed

\ms
Note that the above argument actually gives a general formula for $\th_x(X)$ when
$(S,X)$ is special and condition $(iii)$ of the previous theorem holds. As previously shown in 
\cite[Theorem 3.1]{carkut}, we have
\begin{Cor} Let $G$ be a complex semi-simple algebraic group without any $G_2$ factors.
Then for any $T$-stable subvariety $X$ in $G/P$, we have 
$$\th_x(X)=\bigoplus_{C\in E(X,x)} T_x(C).$$
Consequently, the moment graph $\M(X)$ has exactly $\dim \th_x(X)$ edges at each
vertex $x$.
\end{Cor}

\section{Finite group actions on the moment graph}

Again, suppose $(S,X)$ is special. We will 
next introduce the notion of an action of a finite group $\W$ on $\M(X)$,
and then we will show that an action of $\W$ on $\M(X)$
induces a pair of actions on the $S$-equivariant \coh\ $H^*_S(X^S)$
which leave the image of $H^*_S(X)$ stable as long as $X$ has
vanishing odd \coh.

Let $C\in E(X)$ denote an edge of $\M(X)$, and let  $C^S=\{x,y\}$. 
Suppose $C$ has weights in $X(S)$ given by
$\chi_C$ on $T_x(C)$ and $(\chi_C)^{-1}$ on $T_y(C)$. Taking differentials at the identity,
we obtain an unordered  pair $\{\pm \A\} \subset \s^*$ which determine 
the $\s$-module structure on $T_x(C)$ and $T_y(C)$ respectively. From now on, 
each edge $C$ of $\M(X)$ will be labelled by a pair in this manner. 
We now define a group action on the labelled moment graph. 
\begin{Def} A (finite) group $\W$ is said to {\em act} on the 
moment graph $\M(X)$ of $(S,X)$ if $\W$ acts on  $X^S$, $E(X)$,  and  
linearly on $\s^*$ so that 
if the triple $(x,C,\pm \A)$ determines a labelled edge at $x$, then
 the triple $(w\cdot x,w\cdot C, \pm w\cdot \A)$  
 determines a labelled  edge of $\M(X)$ at $w\cdot x$.
\end{Def}

The following example describes the Weyl group action 
on the  moment graph of the flag variety $\B$.
\begin{Ex}[\cite{cp} and Example 2.5] 
We will first describe the moment graph of $(T,\B)$ and then 
the Weyl group action on $\M(\B)$. 
The action of $W$ on $\B^T$  is by left translation. 
By \cite{cp}, every $T$-stable curve $C$ in $\B$ with $x\in C^T$ has the form
$C=\overline{U_\A x\cdot B}$ where $\A$ is a root of $(G,T)$ 
such that $x^{-1}(\A)<0$, and 
$U_\A$ is the unique one dimensional unipotent subgroup of $G$
normalized by $T$ which is associated to $\A$. Note that 
$C^T=\{x, r_\A x\}$. We claim that this $T$-stable curve has weights 
$\pm \A$. This follows from the fact that there exists an isomorphism 
$\phi:\c \to U_\A$ such that $t\phi(z)t^{-1}=\phi(\A(t)z)$.
We omit the details. 
For $w\in W$, define $w\cdot C=\overline{U_{w(\A)} wx\cdot B}$.
Since $r_{w(\A)}wx=wr_{\A}w^{-1}wx=w r_\A x$, it follows that 
$(w\cdot C)^T=\{wx,wy\}$. Thus if $(x,C,\pm \A)$ determines an edge of $\M(\B)$,
then $(wx,w\cdot C, \pm w(\A))$ determines another edge of $\M(X)$.
The upshot is that as an unlabelled graph, $\M(\B)$ coincides with 
the {\em Bruhat graph $\Gamma(W)$ of} $W$:
namely, the  graph with vertex set $W$, where two vertices $x,y$ are joined by an edge 
$[x,y]$  \iff $yx^{-1}$ is a reflection. 
\end{Ex}

\ms
A special torus action $(S,X)$ where $X$ has vanishing odd \coh\ is 
said to be GKM (\cite{GKM}). If $X$ is smooth, then a special action $(S,X)$ 
is automatically GKM. Of course, the action $(T,\B)$ is GKM. Similarly, 
for any parabolic $P\subset G$, say $T\subset B\subset P$, the
action $(T,G/P)$ is also GKM. In the next example, we show that $W$
acts on $\M(G/P)$.

\begin{Ex}[The $W$-action on $\M(G/P$]
To see that $W$ acts on $\M(G/P)$, notice that
the natural map $\pi:G/B \to G/P$ is a 
$T$-\eq\ regular surjective morphism sending $(G/B)^T$ onto
$(G/P)^T$. For $\pi^{-1}((G/P)^T)$ is closed and $T$-stable,
hence, by the Borel Fixed Point Theorem, it contains a point of $(G/B)^T$.
Now $N_G(T)$ acts transitively on $(G/P)^T$ via $n\to nP$ so that
the stabilzer of $P$ is $N_P(T)$. Thus $(G/P)^T$ is identified with 
$N_G(T)/N_P(T)$, which is itelf identified with $W/W_P$ since $T$
is normal in both  normalizers. The $T$-stable curves in $G/P$ 
are described in a manner similar to those in $G/B$, and it follows 
that $W$ acts on $\M(G/P)$.
\end{Ex}

The next example describes another important GKM-action 
for which $W$ also acts on the moment graph.
\begin{Ex}[Regular Semi-simple Hessenberg Varieties \cite{dps}]\label{SSHV}
Suppose $\h$ is a $B$-submodule of $\g$ containing $\b$.
Fix a regular $t\in \t$, and define a $T$-stable subvariety of $\B$ by 
$$X({\h})=\{ gB\mid g^{-1}\cdot t\in \h\},$$
where, $g\cdot t=\text{Ad}(g)(t).$
One calls $X(\h)$ the {\em Hessenberg variety associated to $\h$.}
Hessenberg varieties are smooth $T$-varieties, and $X(\h)^T=\B^T.$ 
Now $W$ acts on $X(\h)^T$ since $\B^T=W$. 
Moreover, I claim $W$ acts on $\M(X(\h))$. 
To see this, suppose $C\in E(X(\h)),x)$, say $C=\overline{U_\A x\cdot B}$, 
where $\A$ is a root of $(G,T)$ such that $x^{-1}(\A)<0$. 
Now $w\cdot C=\overline{U_{w(\A)} wx\cdot B}=w\overline{U_\A x\cdot B}$.
Let $\phi:\c \to U_\A$ be the isomorphism introduced in the previous example,
and fix $r\in \c^*$. Since $C\subset X(\h)$, it follows that $x^{-1}\phi(-r)\cdot t\in \h$.
Applying the same test to $w\cdot C$, we must check that 
$x^{-1}\phi(-r) w^{-1} \cdot t \in \h$.
But $w^{-1}\cdot t$ is a regular element of $\t$, so $w\cdot C$ is an edge of the moment
graph for the Hessenberg variety $X(\h)$ corresponding to $w^{-1}\cdot t$. 
But the moment graph 
of $X(\h)$ is independent of the choice of a regular $t\in \t$, so $w\cdot C$ is an
edge of  $\M(X(\h))$. It follows that $W$ acts on $\M(X(\h))$. \qed
\end{Ex}

We will now prove the main result of this section:  when $(S,X)$ is GKM, 
every action of a finite group $\W$ on  
a moment graph $\M(X)$ induces a pair of actions on $H^*_S(X^S)$ which 
leave the \ec\ $H^*_S(X)$
invariant.  One of these actions, but not both, satisfy the hypotheses of  Theorem 2.2. To simplify the notation, we will use $C(x,y)$ to denote an $S$-stable curve $C$
with $C^S=\{x,y\}$. Since  $H_S^*(X^S)=\c[\s]\otimes H^0(X^S)$ and $H^0(X^S)=\c^{X^S}$, 
$H_S^*(X^S)$ is naturally identified with $(\c[\s])^{X^S},$ 
the algebra of $\c[\s]$-valued functions on $X^S$.
Supposing $\W$ acts on  $\M(X)$,  define the right and left actions  of $\W$
on $H^*_S(X^S)$  respectively by
$$f\cdot w(x)=f(w\cdot x),$$ 
and 
$$w\bullet f (x)=w\cdot (f(w^{-1}\cdot x)).$$
Note that the right action of $\W$ commutes with the $\c[\s]$-module
structure of $H^*_S(X^S)$ while the left action commutes with a twist. 
  
\begin{Th} 
Suppose $(S,X)$ is a GKM-action,
and suppose $\W$ acts on $\M(X)$. Then the image the 
localization map  $i_S^*:H^*_S(X)\to H^*_S(X^S)$ is stable 
under both the left and right actions of $\W$ on $H^*_S(X^S)$.
Consequently, both of the left and right  actions induce actions of 
$\W$ on $H^*_S(X)$ and consequently on $H^*(X)$ via the forgetful map. 
\end{Th}
\proof 
By the fundamental theorem of \cite{GKM}, the restriction map 
$i_{X}^*: H_S^*(X)\to H_S^*(X^S)$ has image
$$\{ f\in (\c[\s])^{X^S} \mid C(x,y)\in E(X) \implies\, f(x) -f(y) \in \alpha \c[\s],
\mathrm{if}\, C(x,y)\, \mathrm{has ~label}\, \pm\alpha \}.$$
Suppose $f\in H^*_S(X^S)$ lies in the image of $i_S^*$,
and consider an edge  $C=C(x,y)$ of $\M(X)$.  Then
$f(x) -f(y) \in \alpha_{C} \c[\s]$. Thus, 
$$f\cdot w (x) - f\cdot w(y)=f(w\cdot x)-f(w\cdot y) \in w(\alpha_C) \c[\s].$$
Hence the image of $i_S^*$ is stable under the right action.
For the left action, the same reasoning gives
$$w\bullet f(x)-w\bullet f(y)=w\cdot (f(w^{-1}\cdot x)) -w\cdot (f(w^{-1}\cdot y))
\in w(w^{-1}(\alpha_C))\c[\s]=\alpha \c[\s].$$
The final assertion is left to the reader. \qed

\ms 
We will now show the left action of $W$ is trivial on $H^*(G/P)$.
\begin{Th} Let $P$ be an arbitrary parabolic in $G$ containing $B$. 
Then the left action of $W$ on $H^*(G/P)$ is trivial.
\end{Th}
\proof We will first prove the triviality in the case $P=B$. 
Suppose $k$ is a positive integer, and let 
$f\in (\c[\s])^{X^S}$ represent an element of $H^k_T(G/B)$. 
Thus  $f(x)$ is homogeneous of degree $k$ for all $x\in W$. Since $W$ 
is generated by reflections, it suffices to show that
for any root $\A$, $\p(r_\A \bullet f) = \p( f)$,
where $\p:H^*_S(G/B)\to H^*(G/B)$ is the forgetful map.
Now
$r_\A \bullet f (x)=r_\A \cdot f(r_\A x)$, while 
\begin{eqnarray*} r_\A \cdot f(r_\A x)(s)&=& f(r_\A x)(r_\A s)\\
&=&f(r_\A x)(s- \A^{\vee}(s)\A)\\
&\equiv&f(r_\A x)(s) +\A(s)\c[\s]\\
&\equiv&f(x)(s) +\A(s)\c[\s],
\end{eqnarray*}
where $\A^{\vee}=\dfrac{2\A}{<\A,\A>}$.
This implies $r_\A \bullet f$ and $f$ have the same image under the forgetful map, 
so the left action is trivial on $H^*(G/B)$. But by the Leray spectral sequence,
the cohomology map $\pi^*:H^*(G/P)\to H^*(G/B)$ is injective,
since the fibre of the projection map $\pi:G/B \to G/P$ has no odd \coh.
 But $\pi^*$ is \eq, so the theorem is established. 
\qed

\begin{Rem} On the other hand, the left action of $W$ on $H^*(X(\h))$,  
where $X(\h)$ is a regular semi-simple Hessenberg variety is not trivial.
In \cite{tym}, Tymoczko asked how to describe this action.
In  \cite{sw}, Shareshian and Wachs conjectured 
a combinatorial formula for the  character of the left action which involved 
quasisymmetric functions and indifference graphs. This conjecture was 
then proved by Brosnan and Chow in \cite{bc} and by Guay-Paquet \cite{gp}. 
\end{Rem}

\begin{Rem}
For $G/B$, the right action is induced by the usual action of $W$ 
on $K/H$ as mentioned above. The left action was defined  
by Knutson in \cite{knu}. The terms right and left actions were introduced there. 
They were also considered by Tymoczko in \cite{tym} 
where they are called the star and dot actions. 
\end{Rem}

\begin{Rem} Another example of a $T$-variety $X$ in $G/B$ for which $X^T=W$
is a regular $T$-orbit closure. Procesi showed that $H^*(X)$ admits a $W$-action
and classified it in \cite{pro}. It isn't clear that $W$ acts on $\M(X)$.
\end{Rem}

Summarizing the information about the right action from Theorem 2.2 and Corollary 2.3, we get
\begin{Prop}  If $(S,X)$ is a GKM-variety  
and a (finite) group $\W$ acts on the moment graph of $(S,X)$,
then $\W$ acts on $H_S^*(X)$ via the right action by $\c[\s]$-module
automorphisms. Therefore the representation $(\W,H^*(X))$ 
obtained via the forgetful map  is equivalent to the 
representation of $\W$ on $H^0(X^S)$ which is completely determined 
as a sum of induced representations by the orbit structure of $\W$ on $X^S$ as in Corollary \ref{COR1}.
\end{Prop}

FInally, we would like to pose the question of describing the action of $W$ on $H^*(G/P)$
for any parabolic $P$ in $G$.

\ms
\noindent
{\bf Acknowledgement.} I would like to thank Mark Goresky for pointing out the 
paper  \cite{Treu} of Treumann,  and Patrick Brosnan and Kiumars Kaveh for their
many useful remarks.

\end{document}